\documentclass[reqno,centertags, 12pt]{amsart}
\usepackage{amsmath,amsthm,amscd,amssymb}
\usepackage{latexsym}
\usepackage{graphicx}

\newcommand{\bbA}{{\mathbb{A}}}
\newcommand{\bbC}{{\mathbb{C}}}
\newcommand{\bbD}{{\mathbb{D}}}

\newcommand{\bbS}{{\mathbb{S}}}

\newcommand{\dott}{\,\cdot\,}

\newcommand{\lb}{\label}
\newcommand{\f}{\frac}

\newcommand{\ol}{\overline}
\newcommand{\ti}{\tilde  }

\newcommand{\bi}{\bibitem}

\newcommand{\beq}{\begin{equation}}
\newcommand{\eeq}{\end{equation}}
\newcommand{\ba}{\begin{align}}
\newcommand{\ea}{\end{align}}
\newcommand{\veps}{\varepsilon}

\allowdisplaybreaks
\numberwithin{equation}{section}

\newtheorem{theorem}{Theorem}[section]

\newtheorem{proposition}[theorem]{Proposition}

\newtheorem{corollary}[theorem]{Corollary}
\theoremstyle{definition}

\newtheorem{example}[theorem]{Example}

\theoremstyle{remark}

\newcommand{\abs}[1]{\lvert#1\rvert}

\newcounter{smalllist}
\newenvironment{SL}{\begin{list}{{\rm\roman{smalllist})}}{%
\setlength{\topsep}{0mm}\setlength{\parsep}{0mm}\setlength{\itemsep}{0mm}%
\setlength{\labelwidth}{2em}\setlength{\leftmargin}{2em}\usecounter{smalllist}%
}}{\end{list}}

\begin{document}

\title[Fine Structure of the Zeros of OP, II]
{Fine Structure of the Zeros of Orthogonal Polynomials, \\
II. OPUC With Competing Exponential Decay}
\author[B.~Simon]{Barry Simon*}

\thanks{$^*$ Mathematics 253-37, California Institute of Technology, Pasadena, CA 91125, USA. 
E-mail: bsimon@caltech.edu. Supported in part by NSF grant DMS-0140592}

\date{November 8, 2004}

\begin{abstract} We present a complete theory of the asymptotics of the zeros of 
OPUC with Verblunsky coefficients $\alpha_n =\sum_{\ell=1}^L C_\ell b_\ell^n 
+ O((b\Delta)^n)$ where $\Delta <1$ and $\abs{b_\ell}=b<1$. 
\end{abstract}

\maketitle

\section{Introduction} \lb{s1} 

This paper is the second in a series \cite{SaffProc,SaffB3} that discusses detailed 
asymptotics of the zeros of orthogonal polynomials with special emphasis on 
distances between nearby zeros. We will focus here on OPUC, orthogonal polynomials  
on the unit circle; see \cite{SzBk,GBk,OPUC1,OPUC2} for background. The polynomials 
are described by the recursion coefficients, $\{\alpha_n\}_{n=0}^\infty$, called 
Verblunsky coefficients, that give the monic OPUC, $\Phi_n(z)$,  by 
\begin{equation} \lb{1.1} 
\Phi_n(z) = z\Phi_n (z) - \bar\alpha_n \Phi_n^*(z)
\end{equation} 
where 
\begin{equation} \lb{1.2} 
\Phi_n^*(z) = z^n \, \ol{\Phi_n (1/\bar z)}  
\end{equation} 

One of the examples discussed in the first paper \cite{SaffProc} is where $0<b<1$ and 
\begin{equation} \lb{1.3} 
\alpha_n = Cb^n + O((b\Delta)^n) 
\end{equation} 
for $0<\Delta <1$. A typical example is shown in Figure~1 where 
\begin{equation} \lb{1.4} 
\alpha_n = (\tfrac12)^{n+1} 
\end{equation} 
and zeros of $\Phi_{22}$ are shown. Figures and numeric zeros, which appear in 
\cite{OPUC1,SaffProc} and here, are computed using Mathematica and code written by M.~Stoiciu. 

Critical aspects of the zeros in this case are: 
\begin{SL} 
\item[(a)] Finitely many zeros outside $\abs{z}=b+O(\log n/n)$ at the Nevai-Totik 
points, that is, solutions of $D(1/\bar z)^{-1}=0$ where $D$ is the Szeg\H{o} 
function; this is due to Nevai-Totik \cite{NT89}. 

\item[(b)] No zeros in $\abs{z}<b-O(\log n/n)$. This is due to 
Barrios-L\'opez-Saff \cite{BLS}. 

\item[(c)] Zeros near $\abs{z}=b$ are asymptotically a distance $2\pi b/n + 
o(1/n)$ from each other except for a single gap at $z=b$, where the nearby 
zeros are $be^{\pm 2\pi i/n} + O(1/n^2)$, and the distance between these 
neighboring zeros is $2(2\pi b/n) + O(1/n^2)$. This is proven in \cite{SaffProc}. 
\end{SL}

\begin{center}  
\begin{figure}[h] 
\includegraphics[scale=.75]{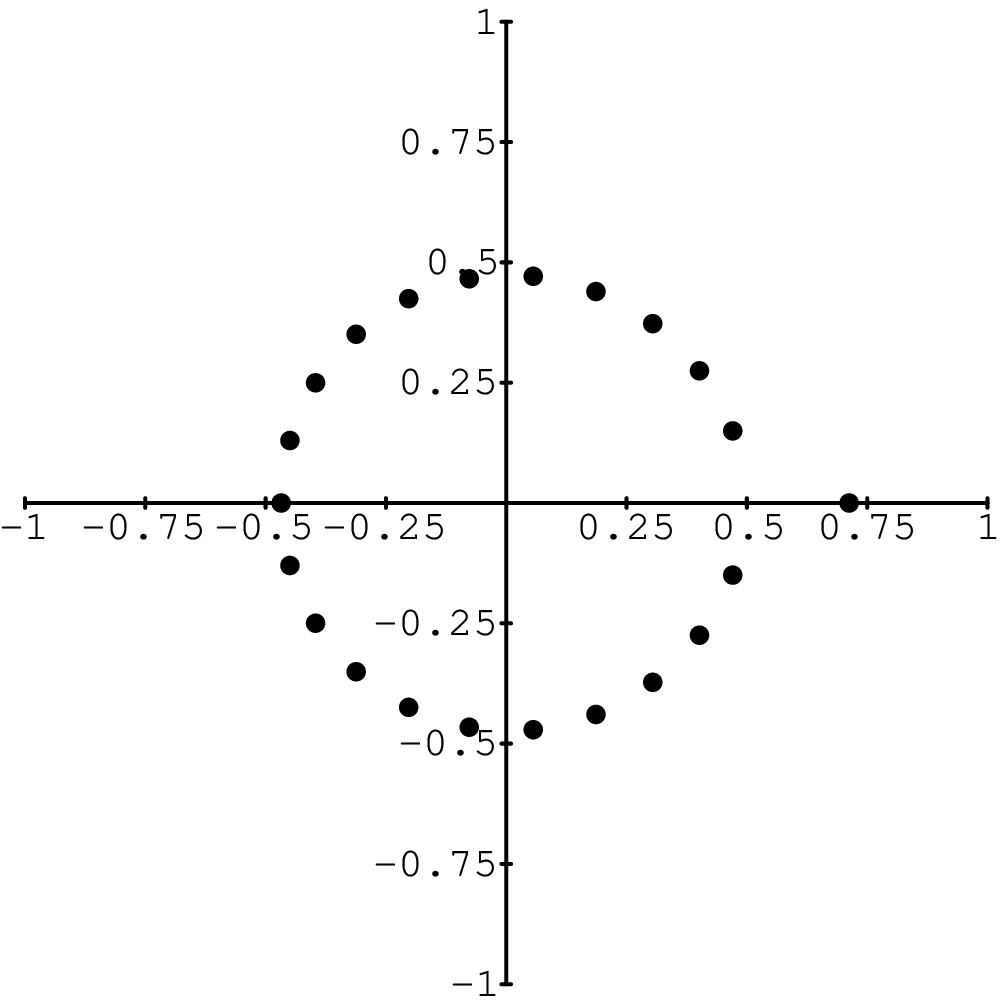}
\caption{}
\end{figure}
\end{center}

In this paper, I will analyze the case of Verblunsky coefficients of the 
form 
\begin{equation} \lb{1.5} 
\alpha_n = \sum_{\ell=1}^L C_\ell b_\ell^n + O((b\Delta)^n)
\end{equation} 
where the $b_\ell$'s are distinct, $C_\ell\neq 0$ for all $\ell$, and 
\begin{equation} \lb{1.6} 
\abs{b_\ell} =b \qquad \ell =1, \dots, L 
\end{equation} 
If the $b_\ell$'s obey $b_\ell^p  =b^p$ for some $p$, then $\alpha_{n+1}/\alpha_n$ 
is periodic of period $p$, and this overlaps examples of BLS \cite{BLS} discussed 
later. 

A typical example is shown in Figure~2 where 
\begin{equation} \lb{1.7} 
\alpha_n = (\tfrac12 )^{n+1} (1+2\cos (\tfrac{\pi}{2} (n+1)))
\end{equation} 
and again, zeros of $\Phi_{22}$ are shown. This has $b_1 = 1/2$, $b_2 = i/2$, 
$b_3 = -i/2$. 

\begin{center}  
\begin{figure} 
\includegraphics[scale=.75]{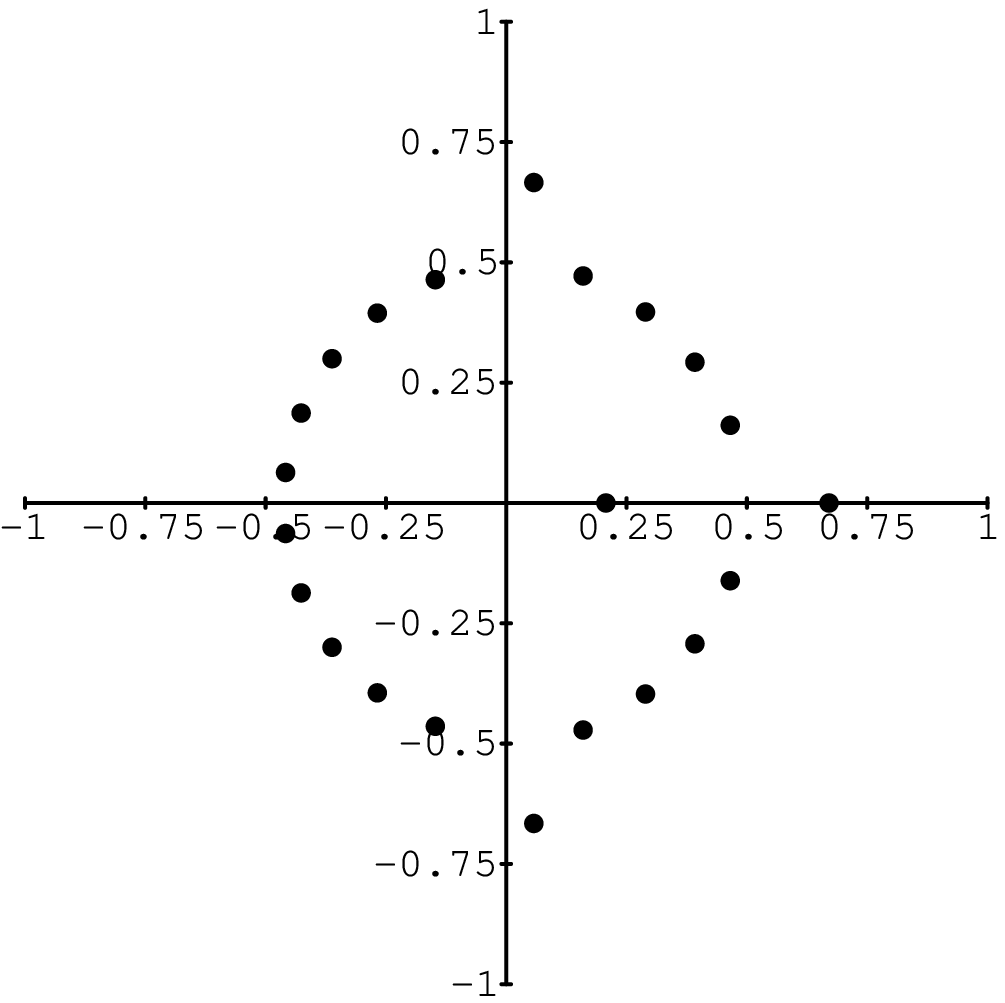}
\caption{}
\end{figure}
\end{center}

Here is what happens to (a)--(c) above: 
\begin{SL} 
\item[(a$^\prime$)] The Nevai-Totik theory still applies. There are three NT zeros 
in the example in Figure~2. 
\item[(b$^\prime$)] There are at most $L-1$ zeros in $\abs{z}<b - O(\log n/n)$ and 
they can be described explicitly. In the case of periodic $\alpha_{n+1}/\alpha_n$, 
the interior zeros for $\Phi_n$ with $n=mp+q$ and $m\to\infty$ have explicitly 
computable limits (limits, not merely accumulation points). For the example in 
Figure~2, $\Phi_n$, $n\equiv 2$ (mod $4$), have a zero approaching $\f12 (\sqrt2 -1) 
=0.20710678118 \dots$. The actual zero in Figure~2 is at $0.20710678374\dots$. 
\item[(c$^\prime$)] There are gaps at each $\bar b_\ell$.  
\end{SL}  

\smallskip
The method used in \cite{SaffProc} to prove (a)--(c) exploits a second-order 
difference equation that $\Phi_n$ obeys, relating $\Phi_{n+1}$ to $\Phi_n$ and 
$\Phi_{n-1}$. That method may extend to the periodic case, $b_\ell^p = b$, but 
will have small divisor problems if $\alpha_n /b^n$ is only almost periodic. 
Instead, this paper will use a different and potentially more powerful and 
illuminating method that views \eqref{1.1} with $\Phi_n^*(z)$ fixed as an 
inhomogeneous first-order difference equation for $\Phi_n(z)$. 

In Section~\ref{s2}, we discuss asymptotics of $\Phi_n$ for $\abs{z}<b$. In 
Section~\ref{s3}, we discuss asymptotics of $\Phi_n$ in the critical region 
$b\Delta_1 < \abs{z} < b\Delta_1^{-1}$. In Section~\ref{s4}, we study zeros in 
$\abs{z}<b$, and in Section~\ref{s5}, the zeros near $\abs{z}=b$. Given 
Section~\ref{s2}, Section~\ref{s4} is straightforward. Section~\ref{s5} will 
use the ideas in \cite{SaffProc}. Finally, Section~\ref{s6} makes various remarks 
about the connection to \cite{BLS}. 

Figures~1 and 2 suggest that there might be a connection between the gaps in
the clock and the Nevai-Totik zeros since in these two cases the number of
NT zeros equals the number of gaps, and the zeros are near the gaps.
There is certainly something to this notion.  If a coefficient $C_\ell$ in 
\eqref{1.5} is changed from zero to nonzero, for large $n$, the zero that 
was in the gap at zero value of $C_\ell$ must stop being on the critical 
circle, and so presumably turns into an NT zero (this is an expectation, 
not a proof, since a proof involves controlling an interchange of limits).  
But the connection is not always there.  Figure~8.3 in \cite{OPUC1} shows an 
example where there is a gap but no associated NT zero.  Moreover, 
the discussion in Section~13 of \cite{SaffProc} makes it clear 
that the number of NT zeros can be arbitrary and is not, in general, 
the same as the number of gaps in the clock.

\medskip 
I thank Percy Deift and Chuck Newman for the hospitality of the Courant Institute 
where some some of this work was done.

\section{Asymptotics in $\abs{z}<b-\veps$} \lb{s2} 

Our main goal here is to give asymptotics of $\varphi_n(z)$ in the region $\abs{z}<b$ 
in case \eqref{1.5} holds. Our methods will also allow us to say something when 
weaker asymptotics (ratio asymptotics) holds and, in particular, to improve 
a result of BLS \cite{BLS}. 

We begin with an analysis of some bounds and rate of convergence of $\Phi_n^*$. 
Estimates similar to these appear in \cite{NT89,MhS1,BLS,SaffProc}. We use 
$q$ rather  than $b$ since sometimes $q=b+\veps$. 

\begin{proposition} \lb{P2.1} Suppose that 
\begin{equation} \lb{2.1} 
\abs{\alpha_n} \leq Cq^n 
\end{equation} 
for some $q\in (0,1)$. Then 
\begin{SL} 
\item[{\rm{(i)}}] With $C_1 =\prod_{j=0}^\infty (1+Cq^j) <\infty$, we have 
\begin{align} 
\abs{z}\leq 1 &\Rightarrow \abs{\Phi_n^*(z)} \leq C_1 \lb{2.2}  \\
\abs{z} \geq 1 & \Rightarrow \abs{\Phi_n(z)} \leq C_1 \abs{z}^n \lb{2.3}
\end{align} 

\item[{\rm{(ii)}}] For $1\leq \abs{z} < q^{-1}$,  
\begin{equation} \lb{2.4} 
\abs{\Phi_n^*(z)} \leq 1+C_1 C\abs{z} (1-q\abs{z})^{-1} 
\end{equation} 

\item[{\rm{(iii)}}] For any $q'>q$, there is $C_{q'}$ with 
\begin{equation} \lb{2.5} 
\abs{\Phi_n(z)} \leq C_{q'} (\max (\abs{z}, q')^n) 
\end{equation} 
for $\abs{z} <1$. 

\item[{\rm{(iv)}}] For $\abs{z} \leq q'$, 
\begin{equation} \lb{2.6} 
\abs{\Phi_n^* (z) -D(z)^{-1} D(0)} \leq \ti C_{q'} (qq')^n 
\end{equation}  

\item[{\rm{(v)}}] $D(z)^{-1}$ has an analytic continuation to $\{z\mid \abs{z} 
<q^{-1}\}$, and in that region, 
\begin{equation} \lb{2.6a} 
\Phi_n^*(z) \to D(0) D(z)^{-1} 
\end{equation} 

\item[{\rm{(vi)}}] For $\abs{z}>q$, 
\begin{equation} \lb{2.6b} 
\lim_{n\to\infty} \, z^{-n} \Phi_n(z) = D(0)\, \ol{D(1/\bar z)}^{-1} 
\end{equation}
\end{SL} 
\end{proposition} 

\begin{proof} (i) For $\abs{z} =1$, $\abs{\Phi_n(z)} = \abs{\Phi_n^*(z)}$, so 
\eqref{2.2} holds by induction from \eqref{1.1}. \eqref{2.2} for $\abs{z}<1$ 
follows from the maximum principle. \eqref{2.3} then follows from \eqref{1.2}. 

\smallskip
(ii) By \eqref{2.3} and the $^*$ of \eqref{1.1}, 
\begin{align} 
\Phi_{n+1}^*(z) &= \Phi_n^*(z) - \alpha_n z\Phi_n(z) \lb{2.7}  \\
\abs{\Phi_{n+1}^*(z)} &\leq 1 + C_1 \sum_{j=0}^n \, \abs{\alpha_j}\, \abs{z}^{n+1} \notag \\ 
&\leq 1 + CC_1 \abs{z} \sum_{j=0}^\infty q^n \abs{z}^n \notag 
\end{align}
proving \eqref{2.4}. 

\smallskip 
(iii) By \eqref{1.2} and \eqref{2.4}, we have \eqref{2.5} for $q' \leq \abs{z} \leq 1$. 
The result for $\abs{z} <q'$ then follows from the maximum principle. 

\smallskip 
(iv) By \eqref{2.7} and \eqref{2.5}, 
\begin{align*} 
\sum_{m=n}^\infty \, \abs{\Phi_{m+1}^*(z) -\Phi_m^*(z)} 
&\leq \sum_{m=n}^\infty C_q C q^m \abs{z} (q')^m \\
&\leq \ti C_{q'} (qq')^n 
\end{align*} 
Since $\Phi_m^*(z) \to D(0) D(z)^{-1}$, \eqref{2.6} holds. 

\smallskip 
(v) Following \cite{NT89}, we note that \eqref{2.7} implies that if $\abs{z}<q^{-1}$, 
then $\sum_n \abs{\Phi_{n+1}^*(z) - \Phi_n^*(z)}<\infty$, so $\Phi_n^*(z)$ has a 
limit. Since Szeg\H{o}'s theorem holds if $\abs{z}<1$, we conclude that $D(z)^{-1}$ 
has a continuation and \eqref{2.6a} holds. 

\smallskip 
(vi) is immediate from (v) and \eqref{1.2}. 
\end{proof} 

By iterating \eqref{1.1}, we obtain  
\begin{equation} \lb{2.8} 
\Phi_n (z) = z^n -\sum_{j=1}^n \bar\alpha_{n-j} z^{j-1} \Phi_{n-j}^*(z) 
\end{equation} 

We conclude that 

\begin{theorem}\lb{T2.2} Suppose \eqref{1.5} holds. Then, uniformly in each disk 
$\abs{z}<b-\veps$, we have that 
\begin{equation} \lb{2.9} 
\begin{split}
\varphi_n(z) - \biggl[\, &\sum_{\ell=1}^L  \bar C_\ell \bar b_\ell^n 
(z-\bar b_\ell)^{-1}  \biggr] D(z)^{-1} \\
&= O(n(b\Delta)^n) + O(nb^{2n}) + O\biggl( nb^n \biggl( 1-\f{\veps}{b}\biggr)^n\biggr)  
\end{split}
\end{equation} 
In particular, uniformly on the disk, 
\[
\lim_{n\to\infty}\, \abs{b^{-n} \varphi_n -\ti Q_n(z)} =0 
\]
where 
\begin{equation} \lb{2.10} 
\ti Q_n(z) =\biggl[\, \sum_{\ell=1}^L \bar C_\ell \omega_\ell^n (\bar b_\ell -z)^{-1} 
\biggr] D(z)^{-1}
\end{equation} 
and  
\begin{equation} \lb{2.11} 
\omega_\ell = \f{\bar b_\ell}{b}  
\end{equation} 
\end{theorem} 

{\it Remarks.} 1. In \eqref{2.9}, the $b^{2n}$ is actually $b^{(3-\delta)n}$ for any 
$\delta >0$. 

\smallskip 
2. $D(z)\ti Q_n(z) \prod_{\ell=1}^L (\bar b_\ell -z)$ is a polynomial of degree $L-1$ 
with almost periodic coefficients (periodic if $b_\ell^p =1$ for some $p$ and all $\ell$).  

\begin{proof} We begin by noting that 
\begin{align} 
\sum_{j=1}^\infty \bar b_\ell^{n-j} z^{j-1} 
&= \bar b_\ell^{n-1} (1-z\bar b_\ell^{-1})^{-1} \notag \\
&= \bar b_\ell^n (\bar b_\ell -z)^{-1} \lb{2.12} 
\end{align} 
and by an identical calculation, 
\begin{equation} \lb{2.13} 
\sum_{j=n+1}^\infty \bar b_\ell^{n-j} z^{j-1} = z^n (\bar b_\ell -z)^{-1}  
\end{equation} 

It follows by \eqref{2.8} that if $\abs{z} <b-\veps$, then 
\begin{equation} \lb{2.14} 
\begin{split}
\biggl| \Phi_n &(z) -D(z)^{-1} D(0) \sum_{\ell=1}^L \bar c_\ell \bar b_\ell^n 
(z-\bar b_\ell)^{-1} \biggr| \\ 
& \leq C_\veps \abs{z}^n + \sum_{j=1}^n \, 
\abs{z}^{j-1} \biggl| \bar\alpha_{n-j} \Phi_{n-j}^* (z) - 
\sum_{\ell=1}^L \bar C_\ell b_\ell^n D(0) D(z)^{-1}\biggr| 
\end{split}
\end{equation} 

The $\abs{z}^n$ term is $O(b^n (1-\f{\veps}{b})^n)$. By \eqref{1.5} and \eqref{2.5}, 
the term in $\abs{\dott}$ in \eqref{2.14} is bounded by 
\begin{equation} \lb{2.15} 
C[b^{2(n-j)} + (b\Delta)^{n-j}] 
\end{equation} 
where the first term comes from $\abs{[\Phi_n^* -D(0) D(z)^{-1}]\alpha_n}$ and the 
second from $\abs{D(z)^{-1} D(0) (\alpha_n -\sum_{\ell=1}^L C_\ell b_\ell^n)}$. 
Since $\abs{z}<b$, the sum in \eqref{2.14} is thus bounded by $nC' [b^{2n} + \max 
(b\Delta,z)^n]$. It follows that 
\[
D(0)^{-1} \Phi_n(z) - \biggl[\, \sum_{\ell=1}^L \bar C_\ell \bar b_\ell^n (z-\bar b_\ell)^{-1} 
\bigg] D(z) = \text{RHS of \eqref{2.9}} 
\]
\eqref{2.9} then follows from 
\begin{align*} 
D(0)^{-1} \Phi_n(z) &= \varphi_n(z) + O\biggl( b^n \biggl|\, \prod_{j=1}^{n-1} 
(1-\abs{\alpha_j}^2)^{1/2} -D(0)\biggr| \biggr) \\ 
&= \varphi_n(z) + O(b^{2n}) 
\qedhere
\end{align*} 
\end{proof}

From \eqref{2.8}, we also get a result that only depends on ratio asymptotics:

\begin{theorem} \lb{T2.3} Suppose that $\alpha_n$ is a sequence of Verblunsky 
coefficients and $n_j$ a subsequence so that 
\begin{SL} 
\item[{\rm{(i)}}] 
\begin{equation} \lb{2.17} 
\limsup_{n\to\infty}\, \abs{\alpha_n}^{1/n} =b\in (0,1)
\end{equation} 
\item[{\rm{(ii)}}] 
\begin{equation} \lb{2.18} 
\liminf_{j\to\infty}\, \abs{\alpha_{n_j}}^{1/n_j} =b
\end{equation} 
\item[{\rm{(iii)}}] For all $k=0,1,2,\dots$ and suitable $\beta_k\in\bbC$, 
\begin{equation} \lb{2.18a} 
\lim_{j\to\infty} \, \f{\bar\alpha_{n_j-k-1}}{\bar\alpha_{n_j-1}} =\beta_k
\end{equation} 
exists.  
\item[{\rm{(iv)}}] For every $\veps$, there is $C_\veps$ so that for all $j$ and 
$k=1,2,\dots$, we have 
\begin{equation} \lb{2.18b} 
\abs{\alpha_{n_j-k-1}} (b-\veps)^k \leq C_\veps \abs{\alpha_{n_j-1}}
\end{equation} 
\end{SL} 

Then for $\abs{z}<b$, $\sum_{j=0}^\infty \beta_j z^j$ converges absolutely and 
uniformly on each disk $\abs{z}<b-\veps$ and 
\begin{equation} \lb{2.19} 
\lim_{j\to\infty}\, \f{\Phi_{n_j}(z)}{\bar\alpha_{n_j-1}} = 
-D(z)^{-1} \sum_{k=0}^\infty \beta_k z^k
\end{equation} 
\end{theorem} 

\begin{proof} By \eqref{2.18b}, we have 
\begin{equation} \lb{2.20} 
\abs{\beta_k} \leq C_\veps (b-\veps)^{-k}
\end{equation} 
proving that $\sum_{j=0}^\infty \beta_j z^j$ has radius of convergence at least $b$. 

In \eqref{2.8} divided by $\alpha_{n_j-1}$, the summand is bounded by 
\[
\abs{z}^{j-1} (b-\veps)^{-j} \sup_{m,\, \abs{z}\leq b}\, \abs{\Phi_m^*(z)} 
\]
which is summable for $\abs{z} <b-\veps$, so by the dominated convergence 
theorem, \eqref{2.19} holds. 
\end{proof} 

\begin{corollary} \lb{C2.4} Let $\alpha_n$ be a sequence of Verblunsky coefficients 
so that 
\begin{equation} \lb{2.21} 
\lim_{n\to\infty} \, \f{\alpha_n}{\alpha_{n-1}} =b\in (0,1)
\end{equation} 
Then uniformly for $\abs{z}<b-\veps$, 
\begin{equation} \lb{2.22} 
\lim_{n\to\infty}\, \f{\Phi_n(z)}{\bar\alpha_{n-1}} = -D(z)^{-1} 
\biggl( 1-\f{z}{b}\biggr)^{-1}
\end{equation} 
\end{corollary} 

{\it Remarks.} 1. By rotational covariance, if \eqref{2.21} holds for some 
$b\in\bbD$, we can find a rotated problem with the ratio in $(0,1)$, so this 
implies a result whenever ratio asymptotics holds. 

\smallskip 
2. This is related to results of BLS \cite{BLS}; see the discussion in 
Section~\ref{s6}. 

\begin{proof} Clearly, \eqref{2.18a} holds with $\beta_k =b^{-k}$ so we need 
only prove \eqref{2.18b}. For any $\delta$, we have 
\[
\biggl| \f{\alpha_{m-1}b}{\alpha_m}\biggr| \leq C_m^{(\delta)} (1+\delta) 
\]
where $C_m^{(\delta)}=1$ for $m\geq M_\delta$ for some $M_\delta$. It follows that 
\[
\biggl| \f{\alpha_{m-k}b^k}{\alpha_m}\biggr| \leq \biggl[ \, \prod_{m=1}^{M_\delta} 
C_m^{(\delta)}\biggr] (1+\delta)^k 
\]
which implies \eqref{2.18b}. 
\end{proof} 

Similarly, we obtain 

\begin{corollary} \lb{C2.5} Let $\alpha_n$ be a sequence of Verblunsky coefficients, 
$b\in (0,1)$, and $c_1, c_2, \dots, c_p$ a sequence so that 
\begin{equation} \lb{2.23} 
\prod_{j=1}^p c_j =1
\end{equation} 
and 
\begin{equation} \lb{2.24} 
\lim_{m\to\infty} \, \f{\alpha_{mp+\ell}}{\alpha_{mp+\ell-1}} = bc_\ell 
\qquad \ell=1,2,\dots, p
\end{equation} 
Then uniformly for $\abs{z}<b-\veps$, 
\[
\lim_{m\to\infty} \, \f{\Phi_{mp+\ell}(z)}{\bar\alpha_{mp+\ell -1}} 
=-D(z)^{-1} G_\ell (z) 
\]
where 
\begin{equation} \lb{2.25} 
\begin{split}
G_\ell (z)\biggl( 1 - \f{z^\ell}{b^\ell}\biggr)  = 1 + (bc_{\ell-1})^{-1} z 
&+ (bc_{\ell-1})^{-1}  (bc_{\ell-2})^{-1} z^2 \\ 
& + \cdots + \prod_{j=1}^{p-1} 
(bc_{\ell-j})^{-1} z^{p-1}
\end{split}
\end{equation}
\end{corollary} 

One can also say something when $b=1$ if we also have 
\begin{equation} \lb{2.29} 
\lim_{n\to\infty} \, \alpha_n =0 
\end{equation}
A key issue is that it may not be true that $\sum_{n=1}^\infty \abs{\alpha_n}^2 
<\infty$, so $D(z)$ may not exist. 

\begin{theorem} \lb{T2.6} Let $\alpha_n$ be a sequence of Verblunsky coefficients 
and $c_1, c_2, \dots, c_p$ is a sequence so that \eqref{2.23} holds and 
\begin{equation} \lb{2.30} 
\lim_{m\to\infty} \, \f{\alpha_{mp+\ell}}{\alpha_{mp+\ell-1}} = c_\ell 
\qquad \ell=1,2, \dots, p
\end{equation} 
Then, uniformly in $\abs{z}<1-\veps$, 
\begin{equation} \lb{2.31} 
\lim_{m\to\infty} \, \f{\varphi_{mp+\ell}(z)}{\bar\alpha_{mp+\ell-1}\varphi_{mp+\ell}^*(z)} 
= -G_\ell(z)
\end{equation} 
where $G$ is given by \eqref{2.25} with $b=1$. 
\end{theorem} 

\begin{proof} \eqref{2.29} and \eqref{2.7}, together with $\abs{\Phi_n(z)} \leq 
\abs{\Phi_n^*(z)}$ on $\bar\bbD$, implies that on $\bar\bbD$, 
\[
\lim_{n\to\infty} \, \f{\Phi_n^*(z)}{\Phi_{n+1}^*(z)} =1 
\]
Dividing \eqref{2.8} by $\bar\alpha_{mp+\ell-1} \Phi_{mp+\ell}^*$, we obtain 
the result by the same argument that led to Corollary~\ref{C2.5}. 
\end{proof}

\section{Asymptotics in the Critical Region} \lb{s3}  

In this section, we will determine asymptotics of $\Phi_n(z)$ in an annulus about 
$\abs{z}=b$ when \eqref{1.5} holds. The idea will be to view \eqref{1.1} as 
an inhomogeneous equation, so we first look at some solutions with particular 
inhomogeneities. Define for $z\neq \bar b_\ell$ and $n=0,1,2,\dots$, 
\begin{equation} \lb{3.1} 
u_n^{(\ell)} = \bar b_\ell^n (z-\bar b_\ell)^{-1} 
\end{equation} 

\begin{proposition} \lb{P3.1} $u_n^{(\ell)}$ obeys 
\begin{equation} \lb{3.2} 
u_{n+1}^{(\ell)} = zu_n^{(\ell)} - \bar b_\ell^n 
\end{equation} 
for all $z\in\bbC$, $z\neq\bar b_\ell$, and all $n=0,1,2, \dots$. 
\end{proposition} 

\begin{proof} $u_{n+1}^{(\ell)} - zu_n^{(\ell)} = (\bar b_\ell^n -z) u_n^{(\ell)} 
=-\bar b_\ell^n$. 
\end{proof} 

Next, define 
\begin{equation} \lb{3.3} 
R_n(z) =\bar\alpha_n \Phi_n^*(z) - \sum_{\ell=1}^L \bar C_\ell b_\ell^n D(z)^{-1} D(0) 
\end{equation} 
and also define 
\begin{equation} \lb{3.4} 
s_n(z) = \sum_{j=0}^\infty z^{-j-1} R_{n+j}(z) 
\end{equation} 
We have 

\begin{proposition} \lb{P3.2} Let $\alpha_n$ obey \eqref{1.5}. Then there is 
$\Delta_1 <1$, 
\begin{SL} 
\item[{\rm{(i)}}] 
\begin{equation} \lb{3.5} 
\sup_{\abs{z} \leq 1}\, \abs{R_n(z)} \leq C (b\Delta_1)^n 
\end{equation} 

\item[{\rm{(ii)}}] The sum in \eqref{3.4} converges uniformly in 
\begin{equation} \lb{3.6x} 
\bbA = \{z\mid 1 > \abs{z} > b\Delta_1\}  
\end{equation} 
and $s_n(z)$ is analytic there. 

\item[{\rm{(iii)}}] We have in $\bbA$ that 
\begin{equation} \lb{3.6} 
\abs{s_n(z)} \leq C(b\Delta_1)^n (\abs{z}-b\Delta_1)^{-1}
\end{equation} 

\item[{\rm{(iv)}}] $s_n$ obeys 
\begin{equation} \lb{3.7} 
s_{n+1}(z) = zs_n(z) - R_n(z) 
\end{equation} 
\end{SL} 
\end{proposition}

\begin{proof} (i) follows from \eqref{1.5} and \eqref{2.6}. 

\smallskip 
(ii), (iii) Since 
\[
\abs{z^{j-1} R_{n+j}(z)} \leq \abs{z}^{-1} (b\Delta_1)^n 
(\abs{z}^{-1} b\Delta_1)^j 
\]
we have a geometric series which yields (ii) and (iii). 

\smallskip 
(iv) Since the sum converges absolutely, 
\begin{align*} 
s_{n+1}(z) - z s_n(z) &= \sum_{j=1}^\infty z^{-j} R_{n+j}(z) - \sum_{j=0}^\infty 
z^{-j} R_{n+j} \\
&= -R_n (z) 
\qedhere 
\end{align*} 
\end{proof} 

The main result of this section is 

\begin{theorem} \lb{T3.3} Let $\alpha_n$ obey \eqref{1.5}. Then for some $\Delta_1 <1$ 
and $z\in\bbA$ given by \eqref{3.6}, we have that 
\begin{equation} \lb{3.8} 
\Phi_n(z) = s_n(z) + \biggl[ \, \sum_{\ell=1}^L \bar C_\ell \bar b_\ell^n 
(z-\bar b_\ell)^{-1} \biggr] D(0) D(z)^{-1} + z^n \ol{D(1/\bar z)}^{-1}
\end{equation} 
\end{theorem} 

{\it Remarks.} 1. Since $\varphi_n = \kappa_n \Phi_n(z)$ and $\kappa_n = D(0)^{-1} 
(1+O(b^n))$, this also gives us asymptotics for $\varphi_n$. 

\smallskip 
2. Since $\Phi_n$ is analytic in $\bbA$, the poles at $\bar b_\ell$ in the second and 
third terms of \eqref{3.8} must cancel. 

\smallskip
3. In $\bbA$, \eqref{3.6} implies $s_n$ is small compared to both $z^n$ and $b^n$, so 
the asymptotics of $\Phi_n$ comes from the competition between the second and third 
terms in \eqref{3.8}. 

\begin{proof} Let 
\[
Q_n(z) =\Phi_n(z) - s_n(z) - \biggl[ \, \sum_{\ell=1}^L \bar C_\ell \bar b_\ell^n 
(z-\bar b_\ell)^{-1} \biggr] D(0) D(z)^{-1} 
\]
By \eqref{1.1}, \eqref{2.2}, and \eqref{3.7}, we have 
\[
Q_{n+1}(z) = zQ_n(z) 
\]
so 
\[
Q_n(z) =f(z) z^n 
\]
Since $Q_n$ is analytic in $\bbA\backslash \{\bar b_\ell\}_{\ell=0}^L$, 
$f(z)$ is analytic there. 

By \eqref{3.6}, 
\[ 
\lim_{n\to\infty}\, \abs{z}^{-n} \abs{s_n(z)} =0 
\]
in $\bbA$, and if $\abs{z}>b$, $\abs{z}^{-n} \sum_{\ell=1}^L \bar C_\ell 
\bar b_\ell^n (z-\bar b_\ell)^{-1}\to 0$, so for $\abs{z} >b$, 
\begin{align*}
f(z) &=\lim_{n\to\infty} z^{-n} Q_n(z) = \lim_{n\to\infty} z^{-n} \Phi_n(z) \\
&= D(0) \, \ol{D(1/\bar z)}^{-1} 
\end{align*} 
by \eqref{2.6b}. 
\end{proof} 

\section{Zeros  in $\abs{z}<b-\veps$} \lb{s4}  

In this section, we use the asymptotic result from Section~\ref{s2} to analyze zeros 
of $\varphi_n$ in the region where $\abs{z}<b$. We initially focus on the case 
where \eqref{1.5} holds. A key role is played by the polynomials 
\begin{equation} \lb{4.1} 
P_n(z) = \sum_{\ell=1}^L \bar C_\ell \omega_\ell^n \prod_{k\neq \ell} 
(z-\bar b_k) 
\end{equation} 
of degree at most $L-1$. Here $\omega_\ell = \bar b_\ell/b$. 

The $P_n$ are almost periodic in $n$ and, in particular, for any sequence $n_j$, 
there is a subsequence $n_{j(k)}$ so $P_\infty \equiv \lim P_{n_{j(k)}}$ exists 
and is a nonzero polynomial (since $P_\infty /\prod_\ell (\bar b_\ell -z)$ has 
poles at each $\bar b_\ell$). 

\begin{theorem}\lb{T4.1} Let \eqref{1.5} hold. Then for any $\veps >0$, there is an 
$N$ so that for $n\geq N$\!, $\varphi_n(z)$ has at most $L-1$ zeros in $\{z\mid \abs{z} 
<b-\veps\}\equiv \bbS$. 
\end{theorem} 

\begin{proof} If not, we can find a sequence $n(j)\to\infty$ so that $P_{n(j)}(z)$ has 
at least $L$ zeros in $\bar\bbS$. By passing to a further subsequence, we can suppose 
$P_{n(j)}\to P_\infty$ and that the $L$ zeros have limits $z_1, \dots, z_L$ in 
$\bar\bbS$ (maybe not distinct). By Theorem~\ref{T2.2}, 
\begin{equation} \lb{4.2} 
\lim_{j\to\infty}\, \varphi_{n(j)} D(z) \prod_{\ell=1}^L (z-\bar b_\ell) = P_\infty(z)  
\end{equation} 
in a neighborhood of $\bar\bbS$, so by Hurwitz's theorem, $P_\infty$ has $L$ zeros 
(counting multiplicity). Since $P_\infty$ has degree $L-1$ and is not identically zero, 
we have a contradiction. 
\end{proof} 

Using Hurwitz's theorem and \eqref{4.2}, we also have an existence result for zeros: 

\begin{theorem}\lb{T4.2} Let \eqref{1.5} hold and let $\omega_\ell = \bar b_\ell/b$. 
Suppose $n(j)$ is a subsequence so that $\lim \omega_\ell^{n(j)}$ exists, call it 
$\omega_\ell^{(\infty)}$. Let 
\begin{equation} \lb{4.3} 
P_\infty (z) =\sum_{\ell=1}^L \bar C_\ell \omega_\ell^{(\infty)} \prod_{k\neq \ell} 
(z-\bar b_k) 
\end{equation} 
and let $\{w_j\}_{j=1}^J$ be its zeros in $\{z\mid\abs{z}<b\}$. Then for all 
sufficiently small $\delta$ and $j\geq N_\delta$, $\varphi_{n(j)}(z)$ has one zero 
within $\delta$ of each $w_j$ and no other zero in $\{z\mid\abs{z} <b-\delta\}$. 
\end{theorem} 

{\it Remark.} By ``one zero within $\delta$ of $w_j$," we actually mean exactly 
$k$ zeros if some $w_j$ occurs $k$ times in the list of zeros counting multiplicity. 

\smallskip 

Since the right side of \eqref{2.22} is nonvanishing on $\{z\mid\abs{z}<b\}$, we recover 
a result of BLS \cite{BLS} from Corollary~\ref{C2.4}, Theorem~\ref{T2.6}, and Hurwitz's theorem: 

\begin{theorem} \lb{T4.3} Let $\alpha_n$ be a sequence of Verblunsky coefficients so 
that 
\[
\lim_{n\to\infty} \, \f{\alpha_n}{\alpha_{n-1}} = b\in (0,1] 
\]
Then for any $\veps >0$, there is an $N_\veps$ so $\varphi_n(z)$ has no zeros in 
$\{z\mid\abs{z}<b-\veps\}$ if $n\geq N_\veps$. 
\end{theorem} 

Finally, Corollary~\ref{C2.5} and Hurwitz's theorem imply 

\begin{theorem}\lb{T4.4} Let $\alpha_n$ be a sequence of Verblunsky coefficients, 
$b\in (0,1)$, and $c_1,c_2, \dots, c_p$ a sequence so that \eqref{2.23} and \eqref{2.24} 
hold. Let $G_\ell$ be given by \eqref{2.22}, let $W_\ell = G_\ell (1-z^\ell/b^\ell) =$  
RHS of \eqref{2.22}, and let $\{w_j^{(\ell)}\}_{j=1}^{N_\ell}$ be the zeros of 
$W_\ell$ in $\{z\mid\abs{z}<b\}$. Then for any 
sufficiently small $\delta$, there is an $N$ so for $mp+\ell \geq N$\!, we have 
that the only zeros of $\varphi_{mp+\ell}$ in $\{z\mid\abs{z}<b-\delta\}$ are one 
each within $\delta$ of each $w_j^{(\ell)}$. 
\end{theorem} 

{\it Remark.} As we will explain in Section~\ref{s6}, that the only possible limit 
points of zeros are the $w_j^{(\ell)}$ is a result of BLS \cite{BLS}, but they 
do not prove there actually are zeros there. 

\begin{example} \lb{E4.5} Let $\alpha_n$ be given by \eqref{1.7}. We have $b=\f12$, 
$p=4$, and $c_1=-1$, $c_2=-1$, $c_3=3$, $c_4=\f13$. Thus 
\[
W_2(z) = 1-2z - 12z^2 - 8z^3 
\]
which has zeros at $-\f12$ and at $\f12 (-1 \pm \sqrt 2)$. Only $(\sqrt 2-1)/2$ is in 
$\{z\mid\abs{z}<\f12\}$. The comparision of the limit and the zeros of $\Phi_{22}$ 
appears in Section~\ref{s1} just after Figure~2. It is not coincidental that $W_2$ 
has a zero at $z=-\f12$. In this case, the second term in \eqref{3.8} is, for 
$n\equiv 2$ (mod $4$), $C(\f12)^n W_2(z)/(z^4 -\f{1}{16})$ with poles only at 
$\f12, \pm \f12 i$. The potential pole at $z=-\f12$ has to be cancelled by a 
zero in $W_2$. 
\qed
\end{example} 

As in \cite{SaffProc}, one can analyze how close the zeros of $\varphi_n$ are to 
the points $w_j^{(\ell)}$. In general, they are exponentially close. If the  
$w_j^{(\ell)}$ are in the annulus where \eqref{3.8} holds, one can write down 
the leading asymptotic exactly. For example, if $w_j^{(\ell)}$ is a $k$-fold zero 
and $D(1/\bar w_j(z))\neq 0$, then the zeros have a clock 
structure as in Theorem~4.5 of \cite{SaffProc}. 

By using Theorem~\ref{T2.6}, we see that Theorem~\ref{T4.4} extends to the case 
$b=1$ if $\lim \alpha_n =0$. In particular, if $\lim \alpha_n =0$ and 
$\lim_{n\to\infty} \f{\alpha_{n+1}}{\alpha_n} =1$ (e.g., $\alpha_n = 
(n+2)^{-\beta}$), then there are no zeros of $\varphi_n$ in $\{z\mid\abs{z} 
<1-\delta\}$ for $n$ large.

\section{Zeros  in the Critical Region} \lb{s5}   

Given Theorem~\ref{T3.3} and the estimate \eqref{3.6}, the analysis of zeros of 
$\varphi_n$ in the region $\{z\mid b\Delta_1 < \abs{z} < b\Delta_1^{-1}\}$ is 
identical to the analysis in \cite{SaffProc} of the zeros in case $L=1$. The gap  
in that case if $z=b$ comes from the analysis of $\ol{D(1/\bar z)}\,\varphi_n(z)$ 
which has zeros with no gap. The gap in zeros of $\varphi_n(z)$ comes from the 
fact that $\ol{D(1/\bar z)}$ has zeros at $z=b$. In our case, when \eqref{1.5} 
holds, $\ol{D(1/\bar z)}$ has a zero at each $\bar b_\ell$, so there are gaps 
at all those points. 

The following extends Theorem~4.3 of \cite{SaffProc} and has the same proof: 

\begin{theorem}\lb{T5.1} Let $\alpha_n$ be a sequence of Verblunsky coefficients 
obeying \eqref{1.5}. Then for some $\delta$, all the zeros $\{z_j^{(n)}\}_{j=1}^{N_n}$ 
of $\varphi_n(z)$ with $\abs{\abs{z}-b}<\delta$ obey 
\begin{SL} 
\item[{\rm{(1)}}] 
\begin{equation} \lb{5.1} 
\sup_j\, \abs{\abs{z_j^{(n)}} -b} = O\biggl( \f{\log n}{n} \biggr)
\end{equation} 
\item[{\rm{(2)}}] For $n$ large, the $z_j^{(n)}$ can be ordered in increasing 
arguments and 
\begin{equation} \lb{5.2} 
\f{\abs{z_{j+1}^{(n)}}}{\abs{z_j^{(n)}}} = 1 + O\biggl( \f{1}{n\log n}\biggr)
\end{equation} 
\item[{\rm{(3)}}] Let $\{\ti z_j^{(n)}\}_{j=1}^{N_n + L}$ be the sequence of 
$z_j^{(n)}$'s with $L$ points added at $\{\bar b_\ell\}_{\ell=1}^L$ still listed 
in increasing order. 
\end{SL} 
Then 
\begin{equation} \lb{5.3} 
\arg z_{j+1}^{(n)} - \arg z_j^{(n)} = \f{2\pi}{n} + O\biggl( \f{1}{n\log n}\biggr) 
\end{equation} 
for $j=1,2,\dots, N_n +L$ with $\arg z_{N_n +L+1}^{(n)}\equiv 2\pi + 
\arg z_1^{(n)}$. Moreover, if $D(z)^{-1}$ is nonvanishing on $\{z\mid 
\abs{z}=b^{-1}\}$, then $O(1/n\log n)$ in \eqref{5.2} and \eqref{5.3} 
can be replaced by $O(1/n^2)$, and $O(\log n/n)$ in \eqref{5.1} can 
be replaced by $O(1/n)$. 
\end{theorem} 

{\it Remark.} In particular, the zeros nearest $\bar b_\ell$ are $\bar b_\ell 
e^{\pm 2\pi i/n} + O(1/n^2)$ with the difference in the args equal to $4\pi /n 
+ O(1/n^2)$. 

\section{Connection to the Results of Barrios-L\'opez-Saff} \lb{s6} 

In this final section, we want to relate the results of \cite{BLS} to ours. In 
their work, determinants of the following form enter:  
\begin{equation} \lb{6.1} 
\Delta_m(z) = \begin{vmatrix} 
\,\, z+x_1 & zx_2 & 0 & {} & {}  \\
1 & z + x_2 & zx_3 & \ddots  \\
0 & 1 & z + x_3 & \ddots & \ddots \\
{} &  {} & \!\!\ddots & \ddots & xz_m \\
{}  & {} & {} & 1 & z +x_m\, \end{vmatrix} 
\end{equation} 
where we also define $\Delta_0 (z) \equiv 1$. We need the following: 

\begin{proposition} 
\begin{SL} 
\item[{\rm{(i)}}] For $m=2,3, \dots$ 
\begin{equation} \lb{6.2} 
\Delta_m(z) = (z+x_m) \Delta_{m-1} - zx_m \Delta_{m-2}(z)  
\end{equation} 
\item[{\rm{(ii)}}] For $m=1,2, \dots$  
\begin{equation} \lb{6.3} 
\Delta_m(z) = z\Delta_{m-1} + x_1 x_2 \dots x_m  
\end{equation} 
\item[{\rm{(iii)}}] For $m=1,2,\dots$ 
\begin{equation} \lb{6.4} 
\Delta_m (z) = z^m + x_1 z^{m-1} + x_1 x_2 z^{m-2} + \cdots + 
x_1 \dots x_m 
\end{equation} 
\end{SL} 
\end{proposition} 

\begin{proof} (i) \eqref{6.2} comes from expanding $\Delta_m$ in minors along the 
last row. 

\smallskip 
(ii) \eqref{6.2} reads 
\[
\Delta_m(z) = z\Delta_{m-1} (z) + x_m [\Delta_{m-1}(z) - z\Delta_{m-2}] 
\]
which implies \eqref{6.3} inductively once one notes that \eqref{6.3} holds for 
$m=1$ since $\Delta_1 (z) = z+x_1$. 

\smallskip 
(iii) This follows by induction from \eqref{6.3}. 
\end{proof} 

In \cite{BLS}, they consider sequences of Verblunsky coefficients where 
\eqref{2.23} and \eqref{2.24} hold to prove that the only accumulation points 
of zeros of $\varphi_{mp+\ell}$ are given by zeros of a polynomial that has 
the form of \eqref{6.1}. Using \eqref{6.4} and 
\[
(x_1\dots x_m)^{-1} \Delta_m (z) = 1+ x_m^{-1} z + x_m^{-1} 
x_{m-1}^{-1} z^2 + \cdots + (x_1 \dots x_m)^{-1} z^m 
\]
one sees their polynomials are up to a constant, our polynomial $W_\ell$. 
Thus our results extend theirs (in that we prove there are, in fact, 
always limit points).

\bigskip

\end{document}